# RESEARCH ANNOUNCEMENT

# MAPPINGS WITH INTEGRABLE DILATATION IN HIGHER DIMENSIONS

JUAN J. MANFREDI AND ENRIQUE VILLAMOR

ABSTRACT. Let $F \in W_{\text{loc}}^{1,n}(\Omega; \mathbb{R}^n)$ be a mapping with nonnegative Jacobian $J_F(x) = \det DF(x) \geq 0$ for a.e. $x$ in a domain $\Omega \subset \mathbb{R}^n$. The *dilatation* of $F$ is defined (almost everywhere in $\Omega$) by the formula

$$K(x) = \frac{|DF(x)|^n}{J_F(x)},$$

Iwaniec and Šverák [IS] have conjectured that if $p \geq n-1$ and $K \in L_{\text{loc}}^p(\Omega)$ then $F$ must be continuous, discrete and open. Moreover, they have confirmed this conjecture in the two-dimensional case $n = 2$. In this article, we verify it in the higher-dimensional case $n \geq 2$ whenever $p > n-1$.

## 1. INTRODUCTION

Let $F: \Omega \to \mathbb{R}^n$ be a mapping in the Sobolev space $W_{\text{loc}}^{1,n}(\Omega; \mathbb{R}^n)$ with non-negative Jacobian $J_F(x) = \det DF(x) \geq 0$ for a.e. $x$ in $\Omega$. We assume that $\Omega$ is a domain in $\mathbb{R}^n$ and that $n \geq 2$. In connection with recent developments in nonlinear elasticity (see [B1, B2, S, IS] and references therein) it has become important to study analytical conditions on the mapping $F$ that will imply topological properties of $F$.

We say that a mapping $F$ has *finite dilatation* if the ratio

$$K(x) = \frac{|DF(x)|^n}{J_F(x)}$$

satisfies

$$1 \leq K(x) < \infty$$

Received by the editors January 4, 1994.

1991 *Mathematics Subject Classification*. Primary 30C65, 35J70 .

*Key words and phrases*. Quasiregular mappings, degenerate elliptic equations, nonlinear elasticity.

First author supported in part by NSF.

for a.e. $x \in \Omega$; that is, if $J_F(x) = 0$ then $DF(x) = 0$, except for a set of measure zero in $\Omega$. Note in particular that if $J_F(x) > 0$ for a.e. $x$, then $F$ has finite dilatation.

It turns out that mappings with finite dilatation are continuous. This was proved by Vodop'yanov and Goldstein [VG] and Šverák [S], using (topological) degree theory. Recently, a new proof has appeared based on the concept of weakly monotone function [M].

If we assume that $K(x) \leq K$ for a.e. $x \in \Omega$, the mapping $F$ is then called *K-quasiregular*. A classical theorem of Rešhetnyak [R] states that $F$ is not only continuous, but also discrete and open.

However, the condition $K \in L_{\text{loc}}^{\infty}$ is too restrictive for the models in elasticity considered in [B1, B2], where typically one looks for minimizers of polyconvex functionals. A simple example of a polyconvex functional is given by

$$(1) \qquad I(F) = \int_{\Omega} \left( |DF(x)|^{\alpha} + \frac{1}{J_F(x)^{\beta}} \right) \, dx.$$

If $F$ satisfies $I(F) < \infty$, and if the condition

$$\frac{n}{\alpha} + \frac{1}{\beta} < \frac{1}{p}$$

holds, it follows from Young's inequality that $K(x) \in L_{\text{loc}}^p$. It is therefore natural to consider whether the boundedness of $K$ is really needed to conclude that $F$ is discrete and open. Suppose that

$$(2) \qquad K \in L_{\text{loc}}^p(\Omega).$$

If $p < n - 1$ there are examples that show that discreteness of $F$ could fail (see [B2], pp. 317–318). For $p \geq n - 1$ and $n = 2$, Iwaniec and Šverák have proved that (2) implies a factorization of Stoilow's type for $F$, from which the discreteness and openness follow immediately. Their proof is based on the existence theorem for quasiconformal mappings in the plane with a given complex dilatation, and it is not clear how to extend it to higher dimensions.

In higher dimensions, $n \geq 2$, Heinonen and Koskela [HK] have proved that if (2) holds for some $p > n - 1$ and $F$ is quasi-light (for each $y \in \mathbb{R}^n$ the components of the set $F^{-1}(y)$ are compact), then $F$ is indeed discrete and open. Without assuming that $F$ is quasi-light they have proved the discreteness and openness for mappings in $W_{\text{loc}}^{1, \frac{4n}{3}}(\Omega, \mathbb{R}^n)$ satisfying (2) for some $p > n - 1$, $n \geq 3$.

Our result states:

**Theorem 1.1.** *Let $F \in W_{loc}^{1,n}(\Omega; \mathbb{R}^n)$ be a mapping of finite dilatation such that for some $p > n - 1$ we have $K \in L_{loc}^p(\Omega)$. Then, either $F$ is constant in $\Omega$ or $F$ is continuous, discrete and open.*

Note that if $F$ has degree one on the boundary, it follows immediately that $F$ is a homeomorphism. Moreover, by a theorem in [HK] it turns out that $F^{-1} \in W^{1,n}$.

In the next two sections we sketch the proof of Theorem 1.1. The details will appear somewhere else [MV].

## 2. First order elliptic systems satisfied by $F$

Our proof is in the spirit of the regularity theory of solutions of degenerate elliptic equations. Denote by $\mathrm{adj}(DF(x))$ the adjugate matrix of $DF(x)$ defined by the relation

$$DF(x) \cdot \mathrm{adj}(DF(x)) = J_F(x).I_n .$$

If $J_F(x) \neq 0$ we have $\mathrm{adj}((DF(x))) = J_F(x)(DF(x))^{-1}$ and, in general, the entries of $\mathrm{adj}(DF(x))$ are homogeneous polynomials of degree $n-1$ with respect to the variables $\frac{\partial F^i}{\partial x_j}$. It is a well-known fact that $\mathrm{adj}(DF)$ satisfies

$$(3) \qquad \mathrm{div}\big(\mathrm{adj}(DF)(V \circ F)\big) = 0$$

in the sense of distributions, where $V$ is a $C^1$ vector field such that $\mathrm{div}\, V = 0$ (see for example [BI]). Setting $V(y) = \overrightarrow{e}_i$ and performing some elementary linear algebra calculations, we obtain that if $u$ is one of the components of $F = (F^1, F^2, \ldots, F^n)$, it satisfies a second-order quasilinear elliptic equation

$$(4) \qquad \mathrm{div}(A(x,\, \nabla u)) = 0\,,$$

where

$$(5) \qquad \frac{1}{c_n K(x)} |\xi|^n \leq A(x,\, \xi) \cdot \xi \leq c_n K^{n-1}(x) |\xi|^n$$

for $\xi \in R^n$. If $K(x) \in L^\infty$, a classical theorem of Serrin [S] gives the Hölder continuity of $F$ right away. To obtain the discreteness and openness, we make the choice $V(y) = \frac{y}{|y|^n}$ in (3) to obtain that $\log |F(x)|$ is a solution, whenever $F(x) \neq 0$, of (4), and is in fact a supersolution everywhere. A weak form of Harnack's inequality for nonnegative supersolutions [BI] implies that

$$\log |F(x)| \in W_{loc}^{1,\, n-\varepsilon} \quad \text{for any } \varepsilon > 0.$$

On applying standard results in nonlinear potential theory, one shows that $F^{-1}(0)$ must have $n$-capacity zero and, thus, Hausdorff dimension zero. It follows that all the connected components of $F^{-1}(0)$ are points. Thus, $F$ is a light mapping (the inverse image of a point is totally disconnected). It now follows from degree theory that a light sense preserving mapping is indeed discrete and open [TY].

When $K$ is not bounded, the inequalities in (5) do not give enough information on (4) to be able to use standard test functions. To proceed in this case, we go back to a full version of (3) when $V$ is not necessarily divergence free.

**Theorem 2.1.** *Let $V$ be a $C^1$-vector field with bounded derivative. Then*

$$(6) \qquad \mathrm{div}((\mathrm{adj}\, DF)\, V \circ F) = [(\mathrm{div}\, V) \circ F]\, J_F$$

*in the sense of distributions.*

*Remark* 1. This formula is quite interesting in itself. It is trivial to check when both $F$ and $V$ are smooth. If we assume a priori that $F$ is quasiregular, it holds for $V \in L^{n/(n-1)}$ as proved in [DS]. For general $F \in W^{1,\, n}$ it can be proved by an approximation argument, and for even more general $F$ it can be found in [MTY].

We now think of (6) as being an elliptic pde and try to prove a supersolution type estimate. The vector field $V(y)$ is chosen to be $|\nabla \Phi(y)|^{n-2} \nabla \Phi(y)$, where $\Phi(y)$ is a $C^2$ "$n$-superharmonic bump function" such that $V(y)$ is $C^1$ and $\Phi(y) \geq \delta > 0$.

Choose $\eta \in C_0^\infty(G)$, $\eta \geq 0$, where $G$ is a relatively compact domain in $\Omega$, and set $\varphi(x) = \eta^n(x)\, \Phi^{1-n}(F(x))$. After standard calculations we get the following estimate of Cacciopoli's type:

$$(7) \qquad \int_G |\nabla(\log \Phi \circ F)(x)|^n \, \eta^n(x) \; \frac{dx}{K(x)} \leq c_n \int_G |\nabla \eta(x)|^n K(x)^{n-1} \, dx.$$

Assume momentarily that we can take $\Phi(y) = \log \frac{1}{|y|}$. A simple calculation using Hölder's inequality and the fact that $K(x) \in L^p$ for some $p > n-1$ gives

$$(8) \qquad \int_G |\nabla(\log \log F)(x)|^{n-1+\epsilon} \eta^{n-1+\epsilon}(x)\, dx < \infty$$

for some $\epsilon > 0$. An argument similar to the case $K \in L^\infty$ gives now that the set $F^{-1}(0)$ has $(n-1+\epsilon)$-capacity zero, and thus, its Hausdorff dimension is less than or equal to $n - (n-1+\epsilon) = 1 - \epsilon < 1$. We conclude that $F$ is light and sense preserving and therefore open and discrete after [TY].

In order to justify the use of $\Phi(y) = \log \frac{1}{|y|}$ we need to approximate $\log \frac{1}{|y|}$ by a sequence $\Phi_a(y)$ of positive, smooth $n$-superharmonic functions and check that the constant $c_n$ in (7) is independent of $a$. We found an explicit approximating sequence by splinning polynomials, and used a computer to do the calculations. This is described in the next section. We are sure that one could also show the existence of such approximation by the shooting method in the theory of nonlinear ode's.

## 3. An $n$-superharmonic bump function

In this section we construct a family of functions $\Phi_a$, each one of them being $n$-superharmonic and smooth, and such that $\Phi_a(y)$ tends to $\log \frac{1}{|y|}$ as $a \to 0$.

Since our arguments are local in nature, without loss of generality we can assume that $F(\Omega) \subset B(0, e^{-e}) = \Omega'$.

**Theorem 3.1.** *For each $0 < a < e^{-e}$, there exists a function $\Phi_a \colon \Omega' \to \mathbb{R}$ with the following properties*:

(i) $\Phi_a \in C^2(\Omega')$;

(ii) $\Phi_a(y) \geq e$, *for every* $y \in \Omega'$;

(iii) $\Phi_a$ *is radial*;

(iv) $\Phi_a'(r) \leq 0$;

(v) $\Phi_a$ *is $n$-superharmonic*;

(vi) $\log \frac{1}{a} \leq \Phi_a(y) \leq \log \frac{1}{a} + \frac{1}{2} + \log 2$, *for every* $|y| \leq e^{-e}$;

(vii) $\Phi_a(y) = \log \frac{1}{|y|}$ *for* $a \leq |y| < e^{-e}$; *and*

(viii) $|\nabla \Phi_a(y)|^{n-2}\, \nabla \Phi_a(y) \ \in \ C^1(\Omega')$.

*Proof.* Here it is!

$$\Phi_a(y) = \begin{cases} \log \frac{1}{|y|}, & \text{if } |y| > a \\ \log \frac{1}{a} - \left(\frac{|y|-a}{a}\right) + \frac{(|y|-a)^2}{2\,a^2}, & \text{if } \frac{a}{2} < |y| < a \\ \log \frac{1}{a} + \log 2 + \frac{1}{2} + (5 - 12\log 2)\, \frac{|y|^2}{a^2} \\ \quad + 4\,(-7 + 12\log 2)\, \frac{|y|^4}{a^4} + 8\,(5 - 8\log 2)\, \frac{|y|^6}{a^6}, & \text{if } |y| < \frac{a}{2}. \end{cases} \qquad \square$$

We find this family by considering polynomials in $|y|$ that agree with $\log \frac{1}{|y|}$ up to order two at $|y| = a$ and look for one that is $n$-superharmonic in $|y| < a$.

Symmetry and smoothness considerations imply that we must consider polynomials in $|y|^2$, making matching at $|y| = a$ not possible.

To overcome this difficulty, we consider the quadratic polynomial $P$ in $|y|$ that agrees with $\log \frac{1}{|y|}$ at $|y| = a$ up to second order for $\frac{a}{2} \leq |y| \leq a$. A quick calculation with the machine reveals that $P$ is $n$-superharmonic. We want a new polynomial $Q$ in $|y|^2$ for $|y| < \frac{a}{2}$ that agrees with $P$ at $|y| = \frac{a}{2}$ up to second order, is $n$-superharmonic and $Q(0) \geq \log \frac{1}{a}$.

It turns out that the value of $Q(0)$ is critical. We now set up the equations that the coefficients of $Q$ must satisfy, considering $Q(0)$ as a parameter, solve symbolically, plot the $n$-Laplacian, and adjust the value of $Q(0)$. We repeat this process until we get a negative plot of the $n$-Laplacian of $Q$.

Once we have an explicit formula for $\Phi_a(y)$, a rigorous proof that it has the required properties is tedious but not very complicated (see [MV]).

*Acknowledgment.* It is a pleasure to thank Juha Heinonen and Pekka Koskela for sharing their work early with us and for many inspiring conversations.

Department of Mathematics, University of Pittsburgh, Pittsburgh, Pennsylvania 15260

*E-mail address*: `manfredi+@pitt.edu`

Department of Mathematics, Florida International University, Miami, Florida 33199

*E-mail address*: `villamor@solix.fiu.edu`